\documentclass[12pt]{article}
\usepackage[a4paper, margin=2.5cm]{geometry}
\usepackage{amsmath,amssymb,amsfonts, color, inputenc}
\usepackage{authblk}
\usepackage{booktabs, url}

\usepackage{latexsym}
\usepackage{amsmath}
\usepackage{amsfonts}
\usepackage{amssymb}
\usepackage{psfrag}
\usepackage{graphicx}

\newtheorem{theorem}{Theorem}

\newtheorem{remark}{Remark}
\newtheorem{definition}{Definition}

\newtheorem{corollary}{Corollary}

\def\bbeta{\boldsymbol{\beta}}

\def\bmu{\boldsymbol{\mu}}
\def\blambda{\boldsymbol{\lambda}}
\def\bDelta{\boldsymbol{\Delta}}
\def\bSigma{\boldsymbol{\Sigma}}

\title{Variable selection for longitudinal survey data}

\author[1]{Laura Dumitrescu }
				
\author[2]{Wei Qian \thanks{Disclaimer. The contents of this report reflect the views of the author and not necessarily the official views or opinions of Statistics Canada.}}
					
\author[3]{J. N. K. Rao \thanks{Research supported by a grant from the Natural Sciences and Engineering Research Council of Canada.}}
\affil[1]{\small{School of Mathematics and Statistics, Victoria University of Wellington, New Zealand}}
\affil[2]{\small{Statistics Canada, Ottawa, Canada}}
\affil[3]{\small{School of Mathematics and Statistics, Carleton University, Canada}}

\begin{document}   
\maketitle

\begin{abstract}
In this article we propose a new variable selection method for analyzing data collected from longitudinal sample surveys. The procedure is based on the survey-weighted quadratic inference function, which was recently introduced as an alternative to the survey-weighted  generalized estimating function. Under the joint model-design framework, we introduce the penalized survey-weighted quadratic inference estimator and obtain sufficient conditions for the existence, weak consistency, sparsity and asymptotic normality. To illustrate the finite sample performance of the model selection procedure, we include a limited simulation study.  
\end{abstract}

\noindent {\em AMS classification}: Primary 62D05; secondary 62F12.

\noindent {\em Keywords}: Complex sampling design; longitudinal data; model selection; oracle property; quadratic inference functions; SCAD; super-population model. 

\section{Introduction}

In economics, social and health sciences, longitudinal sample surveys often exhibit complex sampling design features such as unequal selection probabilities, stratification and clustering of individuals. For data collected from some large-scale surveys, or from surveys which have been linked to administrative data files, to explore relationships between the outcome variables and the covariates, special methods are required for variable selection.   

The problem of variable selection is important since missing important covariates leads to underfitting, inducing estimation bias and poor prediction performance, whereas the inclusion of too many factors in the model generates overfitting, making the model unnecessarily complex and difficult to interpret, as well as producing unstable estimates. Traditional approaches such as stepwise deletion and subset selection are usually used in practice, ignoring the stochastic errors acquired in the variable selection process.

To improve the prediction accuracy and interpretability of regression models with a larger number of covariates, various penalized methods that force some coefficients to zero, have been proposed in the literature. For example, the least absolute shrinkage and selection operator (LASSO) was introduced in \cite{tibshirani96}, but even though it has several attractive properties, for large regression coefficients, the method introduces a significant bias towards zero. Alternative methods have been considered for bias reduction, such as e.g. the adaptive LASSO (ALASSO) (see \cite{zou06}) via a weighted penalty approach, and the smoothly clipped absolute deviation (SCAD) penalty, introduced in \cite{fan-li01}. The latter was shown to satisfy desirable theoretical and empirical properties; for a parameter which is close to zero, it preserves the penalization rate of the LASSO, but, as the absolute value of the parameter  increases, the rate is continuously relaxed. However, the SCAD penalty is non-convex, which generates numerical challenges in obtaining the solution, and, in practice, additional considerations are needed for the selection of the regularization parameter. 

 The traditional approach to analyze sample survey data is to make a design-based inference on the finite population parameters with respect to the distribution induced by the probability sampling design. A penalized method for variable selection, based on an   empirical likelihood approach to model univariate responses from complex sampling surveys, was recently proposed in \cite{zhao-haziza-wu}. Nevertheless, to draw conclusions which are valid beyond the reference population, a stochastic model on the population elements is often needed and in such case, the large sample properties of the estimators are obtained within a joint model-design framework, formally established in \cite{rubin-schiopu05}. The corresponding results specify an average behaviour of estimators which would have been obtained from taking potential samples from all possible finite populations.

The topic of variable selection in models for complex longitudinal sample survey data using penalty functions is rather scarce in the literature and a first procedure, based on the survey-weighted generalized estimating equation (GEE), was introduced in \cite{wang-wang-wang14}. The motivating example was the Canadian National Population Health Survey, where the binary variable of the loss of independence among seniors is modelled using the logistic regression, as a function of eleven other variables from the data set: sex, age, BMI, chronic conditions, smoking status, residence area, education level, income level, living in company, active status and alcohol consumption. To select significant variables and simultaneously estimate coefficients, using the SCAD penalty, the authors used  the survey-weighted GEE, introduced by \cite{rao98} and further discussed in \cite{roberts-ren-rao09} and \cite{carrillo-chen-wu10}.
 
Under the semiparametric marginal modeling approach for longitudinal observations, the correlations between measurements taken at different occasions of the survey are unknown and an alternative to the survey-weighted GEE was proposed in \cite{dumitrescu-qian-rao21}. Its main advantage is that it yields estimators that are more efficient, under misspecification of the working correlation matrix and are as efficient as the GEE counterpart, when the correlation matrix is correctly specified. Moreover, the procedure based on the survey-weighted quadratic inference function (QIF) avoids the additional step of estimating the nuisance correlation parameters and, more importantly, provides an inference function for model diagnostics, as well as for goodness-of-fit tests.

To automatically and simultaneously select variables, we propose the penalized survey-weighted quadratic inference criterion, yielding an estimator with attractive large sample properties. The rates of convergence of the penalized estimator depend on the regularization parameter and, under certain conditions, we show that it retains the oracle property for selecting the correct model: the null components of the estimator are estimated as zero, with probability converging to one, whereas, the nonzero components are estimated as well as the correct submodel is known. 

The paper is organized as follows. Our framework and notations are presented in Section \ref{framework}. In Section \ref{pseudo-QIF-method} we review the survey-weighted GEE and the survey-weighted QIF methods and show how a result in \cite{dumitrescu-qian-rao21} can be used to obtain the limiting distribution of the survey-weighted GEE estimator. The penalized survey-weighted QIF is introduced in Section \ref{penalized-pQIF} and we obtain sufficient conditions for weak consistency, sparsity, as well as the asymptotic normality of its nonzero components. In Section \ref{simulations} we illustrate the finite sample performance of the proposed model selection procedure and obtain numerical results on the estimated coefficients.

\section{The model: assumptions and notations}
\label{framework}

Longitudinal data comprise several observations, made at different time points on a set of individuals or units and recorded measurements consist of a sequence of $n$ size $m$ vectors, denoted as $\{{y}_{i1}, \ldots, y_{im}\},$ $i=1, \ldots, n$. The usual assumption is that there is a correlation within the measurements from each unit $i$ but observations from different units are  independent. Furthermore, corresponding to $y_{ij}$ there is a set of $d$ non-stochastic covariates, denoted as a $d$-dimensional vector ${\bf x}_{ij}$. Hence ${\bf y}_1, \ldots, {\bf y}_n$ is a sample of independent $m$ dimensional random vectors, defined on a probability space $(\Theta, {\cal A}, P_{\beta}),$ $\bbeta \in \Omega$ and the objective is to estimate the parameter $\bbeta$. 

There are several estimating approaches which have been used in the literature and one of the most popular assumes a marginal model  
for the response $y_{ij}$ which depends on the parameter $\bbeta$ through the value $\theta_{ij} = {\bf x}_{ij}^T \bbeta$,  
\begin{equation}
{\rm E}_{\bbeta}(y_{ij}) = \mu(\theta_{ij}):= \mu_{ij}(\bbeta) \mbox{ and } {\rm Var}_{\bbeta}(y_{ij}) = \phi \mu'(\theta_{ij}) := \phi \sigma^2_{ij}(\bbeta), \ \phi \neq 0 \nonumber. \label{marg-model}
\end{equation}
Here, $\mu$ denotes the (canonical) link function, assumed to be a continuously differentiable function, with $\mu'>0$ and $\phi$ is an over-dispersion parameter. In the marginal modelling approach, the true correlation within cluster $i$ is not specified and the GEE  method involves the use of a working correlation matrix instead, denoted as ${\bf R}_i(\boldsymbol{\alpha})$, which depends on a nuisance parameter $\boldsymbol{\alpha}.$ An estimator of $\bbeta$ is taken to be the solution of the GEE, defined as 
\begin{equation}
	{\bf g}_n(\bbeta):=\sum_{i=1}^n \left[\frac{\partial \bmu_i(\bbeta)}{\partial \bbeta^T} \right]^T {\bf A}_i(\bbeta)^{-1/2}{\bf R}_i(\boldsymbol{\alpha})^{-1}{\bf A}_i(\bbeta)^{-1/2}[{\bf y}_i - \bmu_i(\bbeta)]= {\bf 0}, \label{census-GEE}
\end{equation}
where $\bmu_i(\bbeta)=(\mu_{i1}(\bbeta), \ldots, \mu_{im}(\bbeta))^T$ and ${\bf A}_{i}(\bbeta)={\rm diag}\{ \phi \sigma^2_{i1}(\bbeta), \ldots,  \phi \sigma^2_{im}(\bbeta)\}.$ 
Assuming that there exists a ``true'' value of the regression parameter, denoted as $\bbeta_0,$ which is an interior point of $\Omega$, the following identifiability assumption is imposed 
\begin{equation}
    {\rm E}_{\bbeta}[{\bf g}_n(\bbeta)] = {\bf 0} \mbox{ if and only if } \bbeta=\bbeta_0. \label{model-identifiability}
\end{equation}

When the sample is obtained through a complex sampling survey from a finite population $U$ of size $N$, using for example stratification or unequal cluster selection probabilities, the inference needs to account for the sampling features. Namely, from the set of labels $U=\{1, \ldots, N\},$ a subset $s \subset U$ of indices is selected according to a probability sampling design, $P_{\pi}$  such that $P(i \in s) = \pi_i$, $1 =1, \ldots, N$ denotes the first-order inclusion probability. Then, to draw conclusions which are valid beyond the reference population, we assume the stochastic model \eqref{marg-model} on the population elements in $\mathcal{F}=\{{\bf y}_1, \ldots, {\bf y}_N\}$. This leads to a joint model-design inference, where the large sample properties of the estimator are obtained according to the probabilities induced by the model and survey design. The results then specify an average behaviour of the estimator, which would have been obtained from taking potential samples from all possible finite populations. Within the framework of a product probability space of the super-population and the design space (as introduced in \cite{rubin-schiopu05}), an important assumption is that, given the design variables, the sample selection and the model characteristics are independent.
	
	In what follows, the following matrix notations will be used. For a $d \times 1$ vector, ${\blambda},$ we use the notation $\|\blambda\|$ for its Euclidean norm, whereas if ${\bf A}$ is a $d \times d$ matrix, then $\|{\bf A}\| = \sup_{\|\blambda\|=1}\|{\bf A} \blambda\|$ is used for its operator norm. If ${\bf A}$ is symmetric, we denote by $\lambda_{\min}({\bf A})$ and $\lambda_{\max}({\bf A})$ its minimum and maximum eigenvalues, respectively. In addition, for any matrix ${\bf A},$ we have $\|{\bf A}\|=[\lambda_{\max}({\bf A}^T {\bf A})]^{1/2}.$

\section{The survey-weighted quadratic inference function}
\label{pseudo-QIF-method}

A standard assumption in the literature of longitudinal data is that the within-individual correlations are equal and typical choices include 

$(a)\ {\bf R}(\alpha)={\bf I}_m,$ where ${\bf I}_m$ denotes the identity matrix of order $m$,

$(b) \ {\bf R}(\alpha)= \{\rho_{lr}(\alpha)\}_{1 \le l,r \le m},$ where $ \rho_{lr}(\alpha)=\alpha, \ l\neq r, \ 0\le \alpha \le 1,$ $\rho_{ll}(\alpha)=1,$ 

$(c) \ {\bf R}(\alpha)= \{\rho_{lr}(\alpha)\}_{1 \le l,r \le m}, \mbox{ where } \rho_{lr}(\alpha)=\alpha^{|l-r|}, \ 0\le \alpha \le 1,$ 

$(d) \ {\bf R}(\alpha)= \{\rho_{lr}(\alpha)\}_{1 \le l,r \le m}, \mbox{ where } 0 \le \rho_{lr}(\alpha) \le 1,$ $l,r=1, \ldots, m.$

The above forms correspond to $(a)$ {\em independence}, $(b)$ {\em exchangeable}, $(c)$ {\em first-order autoregressive} and $(d)$ {\em unspecified} correlation structure, respectively. 
\begin{remark}
{\rm	As noted in \cite{qu-lindsay-li00}, the inverses of the matrices in (a) - (c) can be written as a linear combinations of a small number of simple basis matrices:

	$(a)$ ${\bf M}_1={\bf I}_m,$  
	
	$(b)$ ${\bf M}_1={\bf I}_m$ and ${\bf M}_2=\{{\gamma}_{lr}\}_{1 \le l,r \le m},$ $\gamma_{lr}=1,$ for $1 \le l \neq r \le m$ and $\gamma_{ll}=0,$ $l=1, \ldots, m$,
	
$(c)$ ${\bf M}_1={\bf I}_m,$ ${\bf M}_2=\{{\gamma}_{lr}\}_{1 \le l,r \le m},$ $\gamma_{l,l-1}=1$, $l=2, \ldots, m$ and $\gamma_{l,r}=0,$for $r \neq l-1,$ $1 \le l, r \le m$ and ${\bf M}_3=\{{\gamma}_{lr}\}_{1 \le l,r \le m},$ with $\gamma_{11}=\gamma_{mm}=1$ and 0 elsewhere.
}	\end{remark}
\vspace{3mm}
Hence, we assume that the working covariance matrix ${\bf R}^{-1}(\boldsymbol{\alpha})$ can be written as 
\begin{equation}
{\bf R}^{-1}(\boldsymbol{\alpha})=\sum_{l=1}^L c_l(\boldsymbol{\alpha}) {\bf M}_l, \label{R-class}
\end{equation}
with given ${\bf M}_1, \ldots, {\bf M}_L$. The ``census extended quasi-score'' vector  $\bar{{\bf q}}_N(\bbeta)$ is defined as  
  \begin{equation}
  \bar{{\bf q}}_N(\bbeta)=\frac{1}{N}\sum_{i=1}^N{\bf q}_i(\bbeta)= \frac{1}{N} \begin{pmatrix}
  \displaystyle{\sum_{i=1}^N  \left[\frac{\partial \bmu_i(\bbeta)}{\partial \bbeta^T} \right]^T {\bf A}_i(\bbeta)^{-1/2}{\bf M}_1{\bf A}_i(\bbeta)^{-1/2}[{\bf y}_i - \bmu_i(\bbeta)]}\\
  \vdots\\
  \displaystyle{\sum_{i=1}^N \left[\frac{\partial \bmu_i(\bbeta)}{\partial \bbeta^T} \right]^T {\bf A}_i(\bbeta)^{-1/2}{\bf M}_L{\bf A}_i(\bbeta)^{-1/2}[{\bf y}_i - \bmu_i(\bbeta)]}
  \end{pmatrix}. \nonumber
  \end{equation}
and the survey-weighted extended quasi-score is obtained using the survey design weights $w_i=\pi_i^{-1},$ $i \in s$ 
\begin{equation}
	{\bf q}_n(\bbeta) =\frac{1}{N} \sum_{i \in s} w_i {\bf q}_i(\bbeta). \nonumber
\end{equation}
If the census extended quasi-score satisfies a central limit theorem (CLT) under the model probability (assumption $(N_1)$) and a CLT holds for the survey-weighted QIF, under the design probability (assumption $(N_2)$), the limiting distribution of the latter was obtained in Theorem 2 of \cite{dumitrescu-qian-rao21}, under the joint model-design probability. The notation $\stackrel{{\cal L}} \longrightarrow$ is used for the convergence in distribution of random variables.

\begin{theorem} \label{CLT-ext-score}
		
	Assume that the following conditions hold.
	
	$(N_0) \ \displaystyle{\frac{n}{N} \longrightarrow f, \ \mbox{as } n \to \infty,}$ with $ 0 \le f <1.$
	
	$(N_1) \ N^{1/2}\bar{{\bf q}}_N(\bbeta_0) \stackrel{{\cal L}} \longrightarrow {\cal N}({\bf 0}, \bSigma_{\beta_0}), \ \mbox{as } N \to \infty, \mbox{ under } P_{\beta_0}, \ \mbox{where }\bSigma_{\beta_0}>0.$
	
	$(N_2) \ \mbox{Given the sequence of finite populations, } \mathcal{F}_N=({\bf y}_{1N}, \ldots, {\bf y}_{NN}),$ we have 
	$$n^{1/2}[{\bf q}_n(\bbeta_0) -\bar{{\bf q}}_N(\bbeta_0)] \stackrel{{\cal L}} \longrightarrow {\cal N}({\bf 0}, \bSigma_{d}), \ \mbox{as } n \to \infty,  \ \mbox{under } P_{\pi}, \ \mbox{where }\bSigma_{d}>0 \mbox{ is non-stochastic}.$$
	Then, as $n \to \infty,$ 
\begin{equation}
\label{CLT-QIF}
	(Q) \  \ n^{1/2}{\bf q}_n(\bbeta_0) \stackrel{{\cal L}} \longrightarrow {\cal N}({\bf 0}, \bSigma_{0}), \mbox{ under } P_{\beta_0,\pi}, \ \mbox{where }\bSigma_{0}=\bSigma_{d}+f\bSigma_{\bbeta_0}.
\end{equation}
\end{theorem}	
 The statement of Theorem \ref{CLT-ext-score} is very general and it can be used to obtain the limiting distribution of a sequence of estimators $\bar{\bbeta}_n$, obtained as a solution of the survey-weighted GEE, ${\bf g}_n^{WGEE}(\bbeta) = {\bf 0}$. The pseudo-GEE estimating function 
\begin{equation}
	{\bf g}_n^{WGEE}(\bbeta):=\sum_{i=1}^n w_i\left[\frac{\partial \bmu_i(\bbeta)}{\partial \bbeta^T} \right]^T {\bf A}_i(\bbeta)^{-1/2}{\bf R}_i(\boldsymbol{\alpha})^{-1}{\bf A}_i(\bbeta)^{-1/2}[{\bf y}_i - \bmu_i(\bbeta)], \label{pseudo-GEE}
\end{equation}
was first proposed in \cite{rao98} and \cite{roberts-ren-rao09} to analyze longitudinal survey data. 

Sufficient conditions for the pseudo-GEE estimator $\bar{\bbeta}_n$ to be weakly consistent, under the joint model and design probability, were given in \cite{carrillo-chen-wu10}. Assuming that $(Q)$ holds for ${\bf g}_n^{WGEE}$, we next show a CLT result for this estimator. 
Denote $\mathcal{B}_{n} = \{\bbeta \in \Omega: \ \  n^{1/2}\|\bbeta - \bbeta_0\| \le r\}, \ \mbox{with }r>0$ and let  
$$ \nabla{\bf g}_n^{WGEE}(\bbeta) = -\frac{\partial {\bf g}_n^{WGEE}(\bbeta)}{\partial \bbeta^T},$$
\begin{equation}
  \bar{{\bf g}}_N(\bbeta)= \frac{1}{N} \sum_{i=1}^N  \left[\frac{\partial \bmu_i(\bbeta)}{\partial \bbeta^T} \right]^T {\bf A}_i(\bbeta)^{-1/2}{\bf R}_i(\boldsymbol{\alpha})^{-1}{\bf A}_i(\bbeta)^{-1/2}[{\bf y}_i - \bmu_i(\bbeta)]. \nonumber
  \end{equation}

\begin{corollary}
Assume that $(N_0)$ holds, together with the following conditions:
	
	$(G_1) \ N^{1/2}\bar{{\bf g}}_N(\bbeta_0) \stackrel{{\cal L}} \longrightarrow {\cal N}({\bf 0}, \bSigma_{\beta_0}^{g}), \ \mbox{as } N \to \infty, \mbox{ under } P_{\beta_0}, \ \mbox{where }\bSigma_{\beta_0}^{g}>0.$
	
	$(G_2) \ \mbox{Given the sequence of finite populations, } \mathcal{F}_N=({\bf y}_{1N}, \ldots, {\bf y}_{NN}),$ we have 
	$$n^{1/2}[{\bf g}_n^{WGEE}(\bbeta_0) -\bar{{\bf g}}_N(\bbeta_0)] \stackrel{{\cal L}} \longrightarrow {\cal N}({\bf 0}, \bSigma_{d}^{g}), \ \mbox{as } n \to \infty,  \ \mbox{under } P_{\pi}, \ \mbox{where }\bSigma_{d}^{g}>0 \mbox{ is non-stochastic}.$$

  $ (G_3) \ \mbox{there exists an invertible non-stochastic } d \times d \mbox{ matrix } {\bf D}_0, \mbox{ such that, as } n \to \infty$, we have 
$$ \sup_{\bbeta \in {\cal B}_n}\|\nabla{\bf g}_n^{WGEE}(\bbeta) - {\bf D}_0\| \stackrel{P_{\bbeta_0, \pi}} \longrightarrow 0.$$ 

Then, as $n \to \infty,$
\begin{equation}
n^{1/2}(\bar{\bbeta}_n - \bbeta_0) \stackrel{{\cal L}} \longrightarrow {\cal N}({\bf 0}, {\bf D}_0^{-1}{\bSigma}_0^{g}{\bf D}_0^{-1}), \ \mbox{under } P_{\bbeta_0,\pi}, \mbox{ where } \bSigma_0^{g} = \bSigma_{d}^{g}+f\bSigma_{\beta_0}^{g}. \label{beta-hat-GEE-CLT}
\end{equation}
\end{corollary}

{\bf Proof}. The Mean Value theorem applied to ${\bf g}_n^{WGEE},$ on the set $\{{\bf g}_n^{WGEE}(\bar{\boldsymbol{\beta}}_n)=0, \bar{\boldsymbol{\beta}}_n \in \mathcal{B}_n\}$ gives  
$${\bf g}_n^{WGEE}(\bar{\boldsymbol{\beta}}_n) ={\bf g}_n^{WGEE}(\boldsymbol{\beta}_0) - \nabla{\bf g}_n^{WGEE}({\boldsymbol{\beta}}_n^*) (\bar{\boldsymbol{\beta}}_n - \boldsymbol{\beta}_0),$$
where ${\boldsymbol{\beta}}_n^* \in \mathcal{B}_n$ and we write 
\begin{eqnarray}
n^{1/2}{\bf g}_n^{WGEE}(\bbeta_0) &=& n^{1/2} {\bf D}_0^{1/2}[ {\bf D}_0^{-1/2}\nabla{\bf g}_n^{WGEE}({\boldsymbol{\beta}}_n^*) {\bf D}_0^{-1/2} -{\bf I}_{d}]{\bf D}_0^{1/2}(\bar{\boldsymbol{\beta}}_n - \boldsymbol{\beta}_0) \nonumber\\
	&+&  n^{1/2}{\bf D}_0(\bar{\boldsymbol{\beta}}_n - \boldsymbol{\beta}_0).\label{CLT-pGEE}
\end{eqnarray} 
By $(G_3),$ the first term in \eqref{CLT-pGEE} is $o_{P_{\bbeta_0 \pi}}(1)$ so that the asymptotic distribution of $n^{1/2}(\bar{\boldsymbol{\beta}}_n - \boldsymbol{\beta}_0)$ is equal to the asymptotic distribution of 
$n^{1/2}{\bf D}_0^{-1} {\bf g}_n^{WGEE}(\bbeta_0)$ and an application of Theorem \ref{CLT-ext-score} concludes the proof. \hfill $\Box$

We now turn to the extended quasi-score vector. Let 
\begin{equation}
	{\bf C}_n(\bbeta) = \frac{1}{N} \sum_{i \in s} w_i {\bf q}_i(\bbeta) {\bf q}_i(\bbeta)^T, \nonumber 
\end{equation}
and assume that it is $P_{\beta_0,\pi}-a.s.$ invertible and that there exists a constant $K>0$ such that for any $n$, we have $\inf_{\bbeta \in \Omega} \lambda_{\min}[{\bf C}_n(\bbeta)] > K, \ P_{\beta_0,\pi}-a.s.$ 

Furthermore, suppose that the link function $\mu$ is  three times continuously differentiable on a neighbourhood of $\bbeta_0,$  $\mathcal{U}_{\delta} = \{\bbeta \in \Omega; \  \|\bbeta - \bbeta_0\| < \delta\}, \ \mbox{with }\delta>0$ so that ${\bf q}_n$ is twice continuously differentiable on $\mathcal{U}_{\delta},$ $P_{\bbeta_0, \pi}$- a.s. Let 
$\displaystyle{ \mathcal{D}_n(\bbeta) = \frac{\partial {\bf q}_n(\bbeta)}{\partial \bbeta^T}},$
whose $k$-th column is given by ${\bf d}_n^{(k)}(\bbeta)$ and, for each $k=1, \ldots, d,$ we denote $\displaystyle{{\bf G}_n^{(k)}(\bbeta) = \frac{\partial {\bf C}_n(\bbeta)}{\partial \beta_k}}$, with uniformly continuous entries on $\mathcal{U}_{\delta}$. 

 The survey-weighted quadratic inference function is defined as 
$$ Q_n(\bbeta)=n{\bf q}_n(\bbeta)^T{\bf C}_n(\bbeta)^{-1}  {\bf q}_n(\bbeta) $$
and its first and second order partial derivatives, with $1 \le k,l \le d$, satisfy
\begin{eqnarray*}
&& \frac{\partial n^{-1}Q_n(\bbeta)}{\partial \beta_k} = 2{\bf d}_n^{(k)}(\bbeta)^T{\bf C}_n(\bbeta)^{-1}{\bf q}_n(\bbeta) - {\bf q}_n(\bbeta)^T{\bf C}_n(\bbeta)^{-1}{\bf G}_n^{(k)}(\bbeta){\bf C}_n(\bbeta)^{-1}{\bf q}_n(\bbeta),   \label{1st-deriv-again}\\
&& \frac{\partial^2 n^{-1}Q_n(\bbeta)}{\partial \beta_k \partial \beta_l} = 2{\bf d}_n^{(k)}(\bbeta)^T{\bf C}_n(\bbeta)^{-1}{\bf d}_n^{(l)}(\bbeta) + r_n^{(k,l)}(\bbeta),
\end{eqnarray*}
where 
\begin{eqnarray*}
&& r_n^{(k,l)}(\bbeta) = 2 \left[\frac{\partial {\bf d}_n^{(k)}(\bbeta)} {\partial \beta_l} \right]^T{\bf C}_n(\bbeta)^{-1}{\bf q}_n(\bbeta) - 2{\bf d}_n^{(k)}(\bbeta)^T{\bf C}_n(\bbeta)^{-1}{\bf G}_n^{(l)}(\bbeta){\bf C}_n(\bbeta)^{-1}{\bf q}_n(\bbeta) +  \\
&& - 2 {\bf d}_n^{(l)} (\bbeta){\bf C}_n(\bbeta)^{-1}{\bf G}_n^{(k)}(\bbeta){\bf C}_n(\bbeta)^{-1}{\bf q}_n(\bbeta) + 2 {\bf q}_n(\bbeta) {\bf C}_n(\bbeta)^{-1} {\bf G}_n^{(l)}(\bbeta) {\bf C}_n(\bbeta)^{-1} {\bf G}_n^{(k)}(\bbeta){\bf C}_n(\bbeta)^{-1}{\bf q}_n(\bbeta) \\
&&  - {\bf q}_n(\bbeta)^T{\bf C}_n(\bbeta)^{-1}\frac{\partial {\bf G}_n^{(k)}(\bbeta)}{\partial \beta_l}{\bf C}_n(\bbeta)^{-1}{\bf q}_n(\bbeta).
\end{eqnarray*}

In \cite{dumitrescu-qian-rao21}, the following assumptions were used to obtain the limiting distribution of the survey-weighted QIF estimator 

$(S_1)$ there exists a non-stochastic and positive definite $Ld \times Ld$ matrix, $ {\cal W}_0(\bbeta)$ such that $\sup_{\bbeta \in {\cal U}_{\delta}}\|{\bf C}_n(\bbeta)^{-1} - {\cal W}_0(\bbeta)\| \stackrel{P_{\beta_0,\pi}} \longrightarrow 0, \mbox{ as } n \to \infty;$ 

$(S_2)$ there exists a non-stochastic matrix $ \mathcal{D}_0(\bbeta),$ of size $Ld \times d,$ such that 
$ \sup_{\bbeta \in {\cal U}_{\delta}}\|\mathcal{D}_n(\bbeta) - \mathcal{D}_0(\bbeta)\| \stackrel{P_{\bbeta_0, \pi}} \longrightarrow 0,$ as $n \to \infty;$
 
$(S_3)$ the matrix ${\cal J}_0(\bbeta) := \mathcal{D}_0(\bbeta)^T{\cal W}_0(\bbeta) \mathcal{D}_0(\bbeta)$ is non-singular on ${\cal U}_{\delta}.$

By Theorem \ref{CLT-ext-score} we have 
\begin{equation}
\sqrt{n}{\bf q}_n(\bbeta_0) = O_{P_{\bbeta_0,\pi}}(1) \nonumber 
\end{equation} 
and, due to $(S_1)$ and $(S_2)$, for every $k =1, \ldots, d$ we obtain 
\begin{equation}
\frac{\partial n^{-1/2}Q_n(\bbeta_0)}{\partial \bbeta}=n^{1/2}\mathcal{D}_n(\bbeta_0)^T{\bf C}_n(\bbeta_0)^{-1}{\bf q}_n(\bbeta) + {\cal R}^{1}_n(\bbeta_0), \ \|{\cal R}^{1}_n(\bbeta_0)\|=o_{P_{\bbeta_0\pi}}(1) . \label{convDQn}
\end{equation}
In addition, since $r_n^{(k,l)}(\bbeta_0) = o_{P_{\bbeta_0, \pi}}(1),$ for any $1 \le k, l \le d$ we have 
 \begin{equation}
	\frac{\partial^2 n^{-1}Q_n({\bbeta}_0)}{\partial \bbeta \partial \bbeta^T} = 2\mathcal{D}_n(\bbeta_0)^T{\bf C}_n(\bbeta_0)^{-1}\mathcal{D}_n(\bbeta_0) + {\cal R}_n^2(\bbeta_0), \mbox{ with } \|{\cal R}_n^2(\bbeta_0)\|=o_{P_{\bbeta_0\pi}}(1), \label{eval-second-der}
	\end{equation}

which implies the element-wise convergence in $P_{\bbeta_0\pi}$ of $\displaystyle{\frac{\partial^2 n^{-1}Q_n({\bbeta}_0)}{\partial \bbeta \partial \bbeta^T}}$ to $2\mathcal{J}_0(\bbeta_0),$ as $n \to \infty$.

\begin{definition}

The pseudo-QIF estimator is defined as 
\begin{equation}
\hat{\bbeta}_n = \arg \min_{\bbeta \in \Omega} Q_n(\bbeta). \label{pseudo-QIF}
\end{equation}
\end{definition}

Under the joint randomization framework, the expected value of the survey-weighted extended quasi-score is equal to 
\begin{equation}
{\rm E}_{\bbeta_0, \pi}[{\bf q}_n(\bbeta)] =\frac{1}{N} \sum_{i=1}^N\bDelta_i(\bbeta), \nonumber
\end{equation}
where
\begin{equation}
\bDelta_i(\bbeta)= {\rm E}_{\bbeta_0}[{\bf q}_i(\bbeta)]= \begin{pmatrix}
\displaystyle{ \left[\frac{\partial \bmu_i(\bbeta)}{\partial \bbeta^T} \right]^T {\bf A}_i(\bbeta)^{-1/2}{\bf M}_1{\bf A}_i(\bbeta)^{-1/2}[\bmu_i(\bbeta_0) - \bmu_i(\bbeta)]}\\
\vdots\\
\displaystyle{ \left[\frac{\partial \bmu_i(\bbeta)}{\partial \bbeta^T} \right]^T {\bf A}_i(\bbeta)^{-1/2}{\bf M}_L{\bf A}_i(\bbeta)^{-1/2}[\bmu_i(\bbeta_0) - \bmu_i(\bbeta)]}
\end{pmatrix}.  \nonumber
\end{equation}
In \cite{dumitrescu-qian-rao21}, the weak consistency of the pseudo-QIF estimator, with respect to the joint model-design probability,  was obtained (see their Theorem 1). Assumptions $(A_1)$ and $(A_2)$ guarantee that the objective function approaches, uniformly, a census-type function, which, due $(A_3)$, has a unique minimum at the true value of the parameter.

\begin{theorem}
\label{pseudo-QIF-consist}

	Assume that the following conditions are satisfied. 
	
	$(A_1)$ $\sup_{\bbeta \in \Omega} \|{\bf C}_n(\bbeta)^{-1} - {\bf W}_N(\bbeta) \| \stackrel{P_{\beta_0,\pi}}\longrightarrow  0, \ \mbox{as } n \to \infty,$ for some positive definite $Ld \times Ld$ matrix ${\bf W}_N(\bbeta)$.

	$(A_2)$ $\displaystyle{\sup_{\bbeta \in \Omega} \left\|{\bf q}_n(\bbeta) - \frac{1}{N} \sum_{i=1}^N\bDelta_i(\bbeta) \right\| \stackrel{P_{\bbeta_0,\pi}} \longrightarrow 0}$, as $n \to \infty.$
		
	$(A_3)$ For every $N \ge 1,$ the equation $\frac{1}{N} \sum_{i=1}^N\bDelta_i(\bbeta)={\bf 0}$ has a unique solution, $\bbeta_0,$

	Then, as $n \to \infty,$ we have $$(C) \  \ \hat{\bbeta}_n \stackrel{P_{\beta_0,\pi}} \longrightarrow {\bbeta_0}.$$
\end{theorem}

The limiting distribution of the pseudo-QIF estimator follows from the convergence of the survey-weighted extended quasi-score (as in Theorem \ref{CLT-QIF}), by showing that this estimator is asymptotically linear, as defined in e.g. (3.3) of \cite{newey-mcfadden86}. When terms of the functions of $\bbeta$ are evaluated at $\bbeta_0,$ we suppress $\bbeta_0$ and denote ${\cal D}_0 = {\cal D}_0(\bbeta_0),$ ${\cal J}_0={\cal J}_0(\bbeta_0)$ and ${\cal W}_0={\cal W}_0(\bbeta_0).$ The next result was shown in Theorem 3 of \cite{dumitrescu-qian-rao21}.

\begin{theorem} 
	\label{pseudo-QIF-CLT}
	
Assume that $(C)$ and $(Q)$ are satisfied, as well as $(S_1)$, $(S_2)$ and $(S_3)$.

Then, as $n \to \infty,$
\begin{equation}
(U) \ \ \sqrt{n}(\hat{\bbeta}_n - \bbeta_0) \stackrel{{\cal L}} \longrightarrow {\cal N}({\bf 0}, {\cal J}_0^{-1}[\mathcal{D}_0^T{\cal W}_0 \bSigma_0 {\cal W}_0\mathcal{D}_0] {\cal J}_0^{-1}), \ \mbox{under } P_{\bbeta_0,\pi}. \label{beta-hat-CLT}
\end{equation}
\end{theorem}
In addition to yielding an estimator which is at least as efficient as the one obtained from the survey-weighted GEE, the survey-weighted QIF can be used to construct a pseudolikelihood ratio type statistics for testing composite hypotheses on model parameters, and a statistic for testing the goodness-of-fit of the marginal model. Their limiting distributions are weighted sums of independent chi-squared random variables, each with one degree of freedom.

\section{The penalized survey-weighted quadratic inference function}
\label{penalized-pQIF}

Based on the survey-weighted QIF, we introduce a new approach to variable selection for longitudinal survey data which can incorporate the within cluster correlation, as well as the survey design features, through penalization: 
$$Q_n^P(\bbeta)=Q_n(\bbeta) + n \sum_{k=1}^d p_{\lambda}(|\beta_k|),$$
where $p_{\lambda}$ is a penalty function which depends on a regularization parameter $\lambda$. The main advantage of penalized methods versus other procedures, such as stepwise deletion and subset selection, is that they can automatically and simultaneously select variables, hence, avoiding the corresponding stochastic errors. There are several functions which have been used as penalties, such as the $L^2$ penalty leading to the ridge regression, or the $L^1$ penalty, used in LASSO. Due to its properties, in our simulations, we use the SCAD, introduced in \cite{fan-li01}, which is a nonconcave penalty function on $(0, \infty)$, defined as a quadratic spline function with knots at $\lambda$ and $a \lambda$ (for some constant $a>0$) 
\begin{equation}
p_{\lambda} (|\theta|) = \begin{cases}
\lambda |\theta|                                              & |\theta| \le \lambda \\
- \frac{|\theta|^2 - 2a\lambda |\theta| + \lambda^2}{2(a-1)} & \lambda < |\theta| \le a \lambda \\
 \frac{(a+1)\lambda^2}{2}                                      & |\theta| > a\lambda.
\end{cases}
\end{equation}
This penalty is singular at origin (yielding sparse solutions), bounded (so that, for large coefficients, the resulting estimators are nearly unbiased) and continuous (leading to a stable model selection procedure). The function is continuously differentiable, with
$$ p'_{\lambda}(|\theta|) = \lambda I(|\theta| \le \lambda) + \frac{a \lambda - |\theta|}{a-1}I(\lambda < |\theta| \le a \lambda ), \ a>2.$$

In \cite{fan-li01}, the value $a=3.7$ was shown to give good practical performance for various selection problems and it was shown that results were similar to those obtained by using the generalized cross-validation method. 

\begin{definition}
\label{penal-pQIF}
The penalized pseudo-QIF estimator is defined as 
$$\tilde{\bbeta}_n = \arg \min_{\bbeta \in \Omega}Q_n^P(\bbeta).$$
\end{definition}	
For any $d$-dimensional vector $\bbeta$, we consider the partition $\bbeta=(\bbeta_1^T, \bbeta_2^T)^T$ into subvectors of size $d_1$ and $d-d_1$, respectively, and use the corresponding notation $$\bbeta_0 = (\bbeta_{10}^T, \bbeta_{20}^T)^T,$$ assuming, without loss of generality that $\bbeta_{20}={\bf 0}$. 

In our simulations, to calculate the penalized pseudo-QIF estimator, we use the iterative method, based on a local quadratic approximation, as proposed in \cite{fan-li01}. The method approximates the nonconvex SCAD penalty term by $\displaystyle{p_{\lambda}(|\beta_k^{(t)}|) + \frac{1}{2} \frac{p'_{\lambda}(|\beta_k^{(t)}|)}{|\beta_k^{(t)}|}[(\beta_k^{(t)})^2 - \beta_k^2] },$ $\beta_k^{(t)}\neq 0$, where $\bbeta^{(t)}=(\beta_1^{(t)}, \ldots, \beta_d^{(t)})^T$ is the estimator obtained at step $t$.  If $\beta_k^{(t)}$ is such that $|\beta_k^{(t)}|<0.001$, we set $\beta_k^{(t+1)}=0$ and write $\bbeta^{(t)}=(\bbeta_1^{(t)}, \bbeta_2^{(t)})^T$, where $\beta_k^{(t)} \neq 0,$ for $k=1, \ldots, d_1$ and $\beta_k^{(t)} = 0,$ for $k=d_1+1, \ldots, d$. The algorithm is initialized with the value of the pseudo-QIF estimator $\hat{\bbeta}_n$ and, based on a previous value $\bbeta^{(t)}$, the objective function is approximated by 
\begin{eqnarray*}
Q_n(\bbeta^{(t)}) + \left[\frac{\partial Q_n(\bbeta^{(t)})}{\partial \bbeta_1}\right]^T(\bbeta_1 - \bbeta_1^{(t)}) &+& \frac{1}{2}(\bbeta_1 - \bbeta_1^{(t)})^T\left[\frac{\partial^2 Q_n(\bbeta^{(t)})}{\partial \bbeta_1 \partial \bbeta_1^T}\right]^T(\bbeta_1 - \bbeta_1^{(t)}) \\
&+& \frac{1}{2} n  \bbeta_1^T \boldsymbol{\Gamma}(\bbeta^{(t)})\bbeta_1,
\end{eqnarray*}
where $\bbeta_1$ is a vector with $d_1$ non-zero entries and $\displaystyle{\boldsymbol{\Gamma}(\bbeta^{(t)})={\rm diag} \left\{\frac{p'_{\lambda}(|\beta_1^{(t)}|)}{|\beta_1^{(t)}|}, \ldots, \frac{p'_{\lambda}(|\beta_{d_1}^{(t)}|)}{|\beta_{d_1}^{(t)}|} \right\}}$. Newton-Raphson algorithm gives the minimizer as 
\begin{equation} 
\bbeta_1^{(t+1)} = \bbeta_1^{(t)} - \left[ \frac{\partial^2 Q_n(\bbeta^{(t)})}{\partial \bbeta_1 \partial \bbeta_1^T} + n \boldsymbol{\Gamma}(\bbeta^{(t)})\right]^{-1} \left[\frac{\partial Q_n(\bbeta^{(t)})}{\partial \bbeta_1}+ n\boldsymbol{\Gamma}(\bbeta^{(t)}) \bbeta_1^{(t)}\right]. \label{NRupdate}
\end{equation}

The performance of the variable selection procedure depends on the choice of the tuning parameter and, following \cite{wang-qu09} and \cite{cho-qu13}, in our simulations, we use a Bayesian information criterion, based on the pseudo-QIF and take $\lambda_n$ to be the minimizer of 
$$ WBIC(\lambda)=[Q_n(\tilde{\bbeta}_{\lambda}) + \log(n) {\rm df}(\tilde{\bbeta}_{\lambda})].$$
Here, $\tilde{\bbeta}_{\lambda}$ denotes the penalized pseudo-QIF estimator and ${\rm df}(\tilde{\bbeta}_{\lambda})$ corresponds to its  number of non-zero entries. It can be shown that, with probability tending to 1, this criterion selects the tuning parameter that identifies the true model (we refer to \cite{qian18} for a proof).

\subsection{Asymptotic properties}
\label{asy-prop-sec}

In this section we investigate the asymptotic properties of the penalized pseudo-QIF estimator $\tilde{\bbeta}_n$. Results are formulated for a general penalty function, whose derivative is continuous and we use the techniques in \cite{fan-li01} to obtain the desired properties. The first theorem shows that, under assumption $(P_1)$, with probability converging to 1, the penalized pseudo-QIF estimator $\tilde{\bbeta}_n$ exists within a ball ${\cal B}_n(r)=\{\bbeta \in \Omega; \ n^{1/2}\|\bbeta - \bbeta_0\| \le r\}$, $r>0$ and hence, it is $\sqrt{n}-$ weakly consistent. Assumption $(P_2)$ ensures that the penalty function does not have much more influence than the pseudo-QIF function on the penalized estimator.

\begin{theorem}{\rm (Existence and weak consistency)}
\label{existence-consistency}

Assume that $(Q)$, $(S_1)$, $(S_2)$, $(S_3)$ are satisfied, together with the following conditions.

$(P_1)$ $ n^{1/2}\max \{p'_{\lambda_n}(|\beta_{0k}|), k=1, \ldots d_1 \} =O(1).$

$(P_2)$ $\max \{|p''_{\lambda_n}(|\beta_{0k}|)|, k=1, \ldots d_1 \} =o(1).$

 Then, there exists a sequence $\tilde{\bbeta}_n$ of random variables satisfying 

${\rm (a)}$ $P_{\bbeta_0, \pi}(\tilde{\bbeta}_n \mbox{ is a local minimizer of }Q_n^P(\bbeta) \mbox { on } {\cal B}_n(r)) \stackrel{n \to \infty}\longrightarrow 1$ and

${\rm (b)}$ $\|\tilde{\bbeta}_n - \bbeta_0\| =O_{P_{\bbeta_0, \pi}}(n^{-1/2}) .$
\end{theorem}
{\bf Proof.}  

(a) We show that, for any $\varepsilon>0$, there exist $r>0$ and $n_{\varepsilon,r}$ such that the event
$\displaystyle{E_n = \left\{Q_n^P(\bbeta_0)<  \inf_{\bbeta \in \partial {\cal B}_n(r)} Q_n^P(\bbeta) \right\} }$ has the property 
\begin{equation}
 P_{\bbeta_0, \pi} (E_n) \ge 1-\varepsilon, \ \mbox{for all } n \ge n_{\varepsilon,r}, \label{exist}
\end{equation}
where $\partial {\cal B}_n(r)=\{\bbeta \in \Omega; \ n^{1/2}\|\bbeta - \bbeta_0\|= r\}.$ This implies that, with probability at least $1-\varepsilon$, there is a minimum in the ball ${\cal B}_n(r)$, i.e. there exists a local minimizer $\tilde{\bbeta}_n$ such that $\|\tilde{\bbeta}_n - \bbeta_0\|=O_{P_{\bbeta_0, \pi}}(n^{-1/2})$.

Let $\bbeta \in \partial{\cal B}_n(r)$ be arbitrarily fixed. By Taylor's expansion of $Q_n(\bbeta)$, using $p_{\lambda_n}(0)=0$
\begin{eqnarray}
Q_n^P(\bbeta) - Q_n^P(\bbeta_0) &=& Q_n(\bbeta) - Q_n(\bbeta_0) + n \sum_{k=1}^d[p_{\lambda_n}(|\beta_k|) - p_{\lambda_n}(|\beta_{0k}|)]\nonumber \\
& \ge & \left[\frac{\partial Q_n(\bbeta_0)}{\partial \bbeta}\right]^T(\bbeta-\bbeta_0) + \frac{1}{2}(\bbeta-\bbeta_0)^T\frac{\partial^2 Q_n(\bbeta^*)}{\partial \bbeta \partial \bbeta^T}(\bbeta-\bbeta_0)\nonumber \\
&+& n \sum_{k=1}^{d_1} \left[p'_{\lambda_n}(|\beta_{0k}|) \frac{|\beta_{0k}|}{\beta_{0k}}(\beta_{k}-\beta_{0k}) +\frac{1}{2} p''_{\lambda_n}(|\beta^{**}_{k}|) (\beta_{k}-\beta_{0k})^2 \right] \nonumber\\
&:=& T_1 + T_2 + T_3 \label{diff}
\end{eqnarray}
where $\bbeta^*$ and $\{\beta^{**}_{k}\}_{1 \le k \le d_1}$ are such that $\|\bbeta^* - \bbeta_0\| \le \|\bbeta - \bbeta_0\| = n^{-1/2}r$ and $|\beta^{**}_k - \beta_{0k}| \le |\beta_k - \beta_{0k}| \le n^{-1/2}r$, for every $k =1, \ldots, d_1.$

The Cauchy-Schwarz inequality gives $\displaystyle{T_1 \ge -r \left \| \frac{\partial n^{-1/2}Q_n(\bbeta_0)}{\partial \bbeta} \right\|}.$

We evaluate 
\begin{eqnarray*}
T_2 &=&  \frac{1}{2}(\bbeta-\bbeta_0)^T \left[\frac{\partial^2 Q_n(\bbeta^*)}{\partial \bbeta \partial \bbeta^T} - \frac{\partial^2 Q_n(\bbeta_0)}{\partial \bbeta \partial \bbeta^T}\right](\bbeta-\bbeta_0) + \frac{1}{2}(\bbeta-\bbeta_0)^T \frac{\partial^2 Q_n(\bbeta_0)}{\partial \bbeta \partial \bbeta^T}(\bbeta-\bbeta_0) \\
&=& o_{P_{\bbeta_0, \pi}}(1) + n(\bbeta-\bbeta_0) ^T {\cal I}_0 (\bbeta-\bbeta_0),
\end{eqnarray*}
where by $(S_3)$, ${\cal I}_0={\cal D}_0^T{\cal W}_0{\cal D}_0$ is positive definite. 

From the Cauchy-Schwarz inequality, using $(P_1)$ and $(P_2)$, we obtain
\begin{eqnarray*}
T_3 & \ge & - r\sqrt{d_1} \sqrt{n}\max\{p'_{\lambda_n}(|\beta_{0k}|), \ k=1, \ldots, d_1\}  - \frac{r^2}{2} \max\{ \left|p''_{\lambda_n}(|\beta_{0k}|)\right|, \ k=1, \ldots, d_1\} + o(1) 
 \end{eqnarray*} 
which is dominated by the leading term of $T_2$.
Hence, using $(P_2)$, we have
\begin{eqnarray*}
Q_n^P(\bbeta) - Q_n^P(\bbeta_0) &\ge& -r \left \| \frac{\partial n^{-1/2}Q_n(\bbeta_0)}{\partial \bbeta} \right\| + r^2 \lambda_{\min}({\cal I}_0)  \\
 &-&  r\sqrt{d_1} \sqrt{n}\max\{p'_{\lambda_n}(|\beta_{0k}|), \ k=1, \ldots, d_1\} - o_{P_{\bbeta_0, \pi}}(1)
\end{eqnarray*}
and for $\varepsilon>0$ arbitrarily fixed, 
choosing $r_{\varepsilon}$ such that, for sufficiently large $n$, the right hand side of the above inequality is strictly positive, with probability of at least $1-\varepsilon$, concludes the proof.

Part (b) now follows from  (a).  \hfill $\Box$

\vspace{2mm}

For the SCAD penalty, if $\lambda_n \to 0$, then, with $n$ large enough, $\max \{p'_{\lambda_n}(|\beta_{0k}|), k=1, \ldots d_1 \} =0$ and $\max \{p''_{\lambda_n}(|\beta_{0k}|), k=1, \ldots d_1 \} =0$ so that $(P_1)$ and $(P_2)$ are satisfied. We now show that, with a regularization parameter chosen such that $\lambda_n \to 0$ and $\sqrt{n} \lambda_n \to \infty$, the penalized pseudo-QIF estimator, using the SCAD penalty performs as well as the oracle procedure. Firstly, it identifies the non-zero components correctly, with a probability converging to 1 and secondly, these estimators are as efficient as the estimator obtained if $\bbeta_{20}={\bf 0}$ were known. In the general case, the assumption $(P_3)$, below assures that the penalty function singular at the origin so that the penalized pseudo-QIF estimator possess the sparsity property.

\begin{theorem}{\rm (Sparsity)}
\label{sparsity}

Assume that $(Q)$, $(S_1)$, $(S_2)$, $(S_3)$ hold, together with 
$$ (P_3) \ \liminf_{n \to \infty} \liminf_{\theta \to 0^+} p'_{\lambda_n}(\theta) \lambda_n^{-1} >0.$$
If $\lambda_n \to 0$ and $\sqrt{n} \lambda_n \to \infty$, as $n \to \infty$, then, with probability converging to 1, for any sequence ${\bbeta}_{1n}$ such that $n^{1/2}\|{\bbeta}_{1n}- \bbeta_{10}\| =O_{P_{\bbeta_0, \pi}}$ and $C>0$, we have $${\bf 0} = \arg \min_{n^{1/2}\|\bbeta_2\| \le C } Q_n[({\bbeta}_{1n}^T, \bbeta_2^T)^T]$$.
\end{theorem}

{\bf Proof.} Let ${\bbeta}_{1n}$ be such that $n^{1/2}\|{\bbeta}_{1n} - \bbeta_{10}\|=O_{P_{\bbeta_0, \pi}}(1)$ and let $C>0$ be arbitrary. 

With $k=d_1+1, \ldots, d$, using a Taylor's expansion of the function $\displaystyle{\frac{\partial Q_n^{P}[({\bbeta}_{1n}^T, \bbeta_2^T)^T]}{\partial \beta_k}}$ (as a function of $\bbeta_2$), around $\bbeta_{20}={\bf 0}$, we have
\begin{eqnarray}
\frac{\partial Q_n^{P}[({\bbeta}_{1n}^T, \bbeta_2^T)^T]}{\partial \beta_k} &=& \frac{\partial Q_n[({\bbeta}_{1n}^T, {\bf 0}^T)^T]}{\partial \beta_k} + \sum_{l=d_1+1}^{d} \frac{\partial^2Q_n[({\bbeta}_{1n}^T, \bbeta_2^{*T})^T]}{\partial \beta_k \partial \beta_l}\beta_l + n p'_{\lambda_n}(|\beta_k|) \mbox{sign}(\beta_k), \nonumber\\
&=& \frac{\partial Q_n[(\bbeta_{10}^T, {\bf 0}^T)^T]}{\partial \beta_k} +o_{P_{\bbeta_0 \pi}}(1) + \sum_{l=d_1+1}^{d}\frac{\partial^2Q_n[(\bbeta_{10}^T, {\bf 0}^T)^T]}{\partial \beta_k \partial \beta_l}\beta_l +o_{P_{\bbeta_0 \pi}}(1)\nonumber\\
&+& n p'_{\lambda_n}(|\beta_k|) \mbox{sign}(\beta_k), \label{decomp}
\end{eqnarray}
where $\bbeta_2^*$ is such that $\|\bbeta_2^* \| \le \|\bbeta_2\| \le Cn^{-1/2}$ and we used the continuity of the first and second order partial derivatives of $Q_n(\bbeta)$, with respect to $\beta_k$, $k=d_1+1, \ldots, d$ around $\bbeta_0$. From \eqref{convDQn}, the first term in \eqref{decomp} is $O_{P_{\bbeta_0 \pi}}(n^{1/2}),$ whereas \eqref{eval-second-der}, together with $(S_1)$ and $(S_2)$, implies that the third  term is $O_{P_{\bbeta_0 \pi}}(n^{1/2}).$ 

Since $n^{-1/2}\lambda_n^{-1} \to 0$, we obtain
\begin{eqnarray*} 
\frac{\partial Q_n^{P}[({\bbeta}_{1n}^T, \bbeta_2^T)]}{\partial \beta_k} &=& n p'_{\lambda_n}(|\beta_k|) \mbox{sign}(\beta_k) + O_{P_{\bbeta_0\pi}}(n^{1/2}) + o_{P_{\bbeta_0\pi}}(1)\\
&=&  n \lambda_n   [\lambda_n^{-1}p'_{\lambda_n}(|\beta_k|) \mbox{sign}(\beta_k)  + o_{P_{\bbeta_0\pi}}(1) + o_{P_{\bbeta_0\pi}}(n^{-1} \lambda_n^{-1})] ,
\end{eqnarray*}
which, due to $(P_3)$, for large enough $n$, we have $\displaystyle{\mbox{sign}\left\{\frac{\partial Q_n^{P}[({\bbeta}_{1n}^T, \bbeta_2^T)]}{\partial \beta_k}\right	\} = \mbox{sign}(\beta_k)}$ for any $k =d_1+1, \ldots, d$ and $C>0$. Hence, $\displaystyle{Q_n^{P}[({\bbeta}_{1n}^T, \bbeta_2^T)]}$ has a local minimum within the ball $\{\bbeta_2: \ n^{1/2}\|\bbeta_2\| \le C\}$, at $\bbeta_{2}={\bf 0}$. \hfill $\Box$

\vspace{2mm}

Under the assumptions of Theorem \ref{sparsity}, with probability converging to 1, if the penalized pseudo-QIF estimator $\tilde{\bbeta}_n = (\tilde{\bbeta}_{1n}^T, \tilde{\bbeta}_{2n}^T)^T$ is $\sqrt{n}$-consistent, then it must satisfy $\tilde{\bbeta}_{2n}={\bf 0}$. 

Let ${\cal D}_{10}$ be the $Ld_1 \times d_1$ matrix obtained from ${\cal D}_0\left\{\begin{pmatrix}
{\bbeta}_{10} \\
{\bf 0}
\end{pmatrix}\right\}$ by selecting the rows and columns corresponding to $\bbeta_{10}$, i.e. include the first $d_1$ columns and rows numbered by $\eta d+1, \ldots, \eta d+d_1$, where $\eta=0, \ldots, (L-1)$. Similarly, let ${\cal W}_{10}$ be the $Ld_1 \times Ld_1$ matrix obtained from ${\cal W}_0\left\{\begin{pmatrix}
{\bbeta}_{10} \\
{\bf 0}
\end{pmatrix}\right\}$ by selecting the rows and columns which are numbered as $\eta d+1, \ldots, \eta d+d_1$, where $\eta=0, \ldots, (L-1)$. Finally, let ${\cal I}_{10}$ denote the upper $d_1 \times d_1$ corner of the matrix ${\cal I}_0\left\{\begin{pmatrix}
{\bbeta}_{10} \\
{\bf 0}
\end{pmatrix}\right\}$ and 
\begin{eqnarray*}
{\bf b}_n &=&  \left(p'_{\lambda_n}(|\beta_{01}|)\mbox{sign}(\beta_{01}), \ldots, p'_{\lambda_n}(|\beta_{0d_1}|)\mbox{sign}(\beta_{0d_1}) \right)^T, \\
{\bf B}_{n} &=&  {\rm diag}\{p''_{\lambda_n}(|\beta_{01}|), \ldots, p''_{\lambda_n}(|\beta_{0d_1}|)\}.\\
\end{eqnarray*}
 
Theorem \ref{normality} shows that the limiting distribution of $\tilde{\bbeta}_{1n}$ is normal and that, for certain penalties, including the SCAD, the penalized pseudo-QIF of $\bbeta_{10}$ is asymptotically as efficient as the estimator obtained if $\bbeta_{20}={\bf 0}$ were known.

\begin{theorem}{{\rm (Limiting distribution)}}
\label{normality}

Assume that $(Q)$, $(S_1)$, $(S_2)$, $(S_3)$ and $(P_3)$ hold. If $\lambda_n \to 0$ and $\sqrt{n} \lambda_n \to \infty$, as $n \to \infty$, then, with probability converging to 1, the $\sqrt{n}$-consistent local minimizer of $Q_n^{P}(\bbeta)$, $\tilde{\bbeta}_n=(\tilde{\bbeta}_{1n}^T, {\bf 0}^T)^T$ satisfies
 $$\sqrt{n}[2{\cal I}_{10} + {\bf B}_{n}]\{[\tilde{\bbeta}_{1n} - \bbeta_{10} + [2{\cal I}_{10} + {\bf B}_{n}]^{-1} {\bf b}_n\} \stackrel{{\cal L}} \longrightarrow {\cal N}({\bf 0}, 4{\cal D}_{10}^T {\cal W}_{10}\bSigma_{0}^{(1)}{\cal W}_{10}{\cal D}_{10}),$$ under $P_{\bbeta_0 \pi}$,  where $\bSigma_{0}^{(1)}$ denotes the $Ld_1 \times Ld_1$ matrix obtained from $\bSigma_0=\bSigma_d + f \bSigma_{(\bbeta_{10}^T, {\bf 0}^T)^T}$ by selecting the rows and columns which are numbered as $\eta d+1, \ldots, \eta d+d_1$ and $\eta=0, \ldots, (L-1).$
\end{theorem}
{\bf Proof.} Since the local minimizer of $Q_n^P(\bbeta_{1}, {\bf 0})$, $\tilde{\bbeta}_{1n}$, is a $\sqrt{n}$-consistent estimator of $\bbeta_{10},$ a Taylor's expansion gives
\begin{eqnarray*}
\frac{\partial Q_n^{P}}{\partial \bbeta_1} \left\{\begin{pmatrix}
\tilde{\bbeta}_{1n} \\
{\bf 0}
\end{pmatrix}\right\} &=& \frac{\partial Q_n}{\partial \bbeta_1} \left\{\begin{pmatrix}
{\bbeta}_{10} \\
{\bf 0}
\end{pmatrix} \right\} + \frac{\partial^2 Q_n}{\partial \bbeta_1 \partial \bbeta_1^T} \left\{\begin{pmatrix}
{\bbeta}_{1}^* \\
{\bf 0}
\end{pmatrix} \right\}(\tilde{\bbeta}_{1n}-\bbeta_{10}) +n{\bf b}_n \\
&+& n {\rm diag}\{p''_{\lambda_n}(|\beta_{1}^{**}|), \ldots, p''_{\lambda_n}(|\beta_{d_1}^{**}|)\}(\tilde{\bbeta}_{1n}-\bbeta_{10}).
\end{eqnarray*} 
where $\bbeta_1^*$ and $\beta_{k}^{**}$ are such that $\|\bbeta_1^* - \bbeta_{10}\| \le \|\tilde{\bbeta}_{1n} - \bbeta_{10}\|$ and $|\beta_{k}^{**} - \beta_{0k}| \le |\tilde{\beta}_{nk} - \beta_{0k}|,$ with $k =1 \ldots, d_1$.
After rearranging, we obtain
\begin{eqnarray*}
&&-\sqrt{n}\frac{\partial n^{-1}Q_n}{\partial \bbeta_1} \left\{\begin{pmatrix}
{\bbeta}_{10} \\
{\bf 0}
\end{pmatrix} \right\} = \frac{\partial^2 n^{-1}Q_n}{\partial \bbeta_1 \partial \bbeta_1^T} \left\{\begin{pmatrix}
{\bbeta}_{10} \\
{\bf 0}
\end{pmatrix} \right\}\sqrt{n}(\tilde{\bbeta}_{1n}-\bbeta_{10}) + \sqrt{n}{\bf b}_n + {\bf B}_n \sqrt{n}(\tilde{\bbeta}_{1n}-\bbeta_{10}) \\
&& + \left[\frac{\partial^2 n^{-1}Q_n}{\partial \bbeta_1 \partial \bbeta_1^T}\left\{\begin{pmatrix}
{\bbeta}_{1}^* \\
{\bf 0}
\end{pmatrix} \right\} - \frac{\partial^2 n^{-1}Q_n}{\partial \bbeta_1 \partial \bbeta_1^T}\left\{\begin{pmatrix}
{\bbeta}_{10} \\
{\bf 0}
\end{pmatrix} \right\}  \right]\sqrt{n}(\tilde{\bbeta}_{1n}-\bbeta_{10}) \\
&& + {\rm diag} \{p''_{\lambda_n}(|\beta_{1}^{**}|) - p''_{\lambda_n}(|\beta_{01}|), \ldots, p''_{\lambda_n}(|\beta_{d_1}^{**}|)-p''_{\lambda_n}(|\beta_{0d_1}|)\}\sqrt{n}(\tilde{\bbeta}_{1n}-\bbeta_{10}).
\end{eqnarray*}
 Due to \eqref{eval-second-der}, $\displaystyle{\frac{\partial^2 n^{-1}Q_n}{\partial \bbeta_1 \partial \bbeta_1^T} \left\{\begin{pmatrix}
{\bbeta}_{10} \\
{\bf 0}
\end{pmatrix} \right\} = 2{\cal I}_{10} + o_{P_{\bbeta_0\pi}(1)}}$ which, together with the fact that the entries of $\displaystyle{\frac{\partial^2 n^{-1}Q_n}{\partial \bbeta_1 \partial \bbeta_1^T}\left\{\begin{pmatrix}
{\bbeta}_{1} \\
{\bf 0}
\end{pmatrix}\right\}}$ are continuous at $(\bbeta_{10}^T,{\bf 0}^T)^T$, the asymptotic distribution of 
$$ \sqrt{n} [2{\cal I}_{10}+{\bf B}_n] \{\tilde{\bbeta}_{1n}-\bbeta_{10} + [2{\cal I}_{10}+{\bf B}_n]^{-1}{\bf b}_n\}$$ is equivalent to that of 
$\displaystyle{-\sqrt{n}\frac{\partial n^{-1}Q_n}{\partial \bbeta_1} \left\{\begin{pmatrix}
{\bbeta}_{10} \\
{\bf 0}
\end{pmatrix} \right\}}$. The latter is equivalent to the limiting distribution of $\displaystyle{-2\frac{\partial {\bf q}_n}{\partial \bbeta_1^T}\left\{\begin{pmatrix}
{\bbeta}_{10} \\
{\bf 0}
\end{pmatrix} \right\} \left\{{\bf C}_n\left\{\begin{pmatrix}
{\bbeta}_{10} \\
{\bf 0}
\end{pmatrix} \right\}\right\}^{-1} \sqrt{n}{\bf q}_n\left\{\begin{pmatrix}
{\bbeta}_{10} \\
{\bf 0}
\end{pmatrix} \right\} }$, using \eqref{convDQn}, which, in turn, is given by ${\cal N}({\bf 0}, 4{\cal D}_{10}^T{\cal W}_{10} \bSigma_{0}^{(1)}{\cal W}_{10}{\cal D}_{10})$, under the joint model-design probability, due to Theorem \ref{CLT-ext-score} and assumptions $(S_1)$ and $(S_2)$. \hfill $\Box$

\begin{remark}
{\rm Under the assumptions of Theorem \ref{normality}, if $\|{\bf b}_n\| \to 0$ and $\|{\bf B}_n\| \to 0$ then 
$$\sqrt{n}(\tilde{\bbeta}_{1n} - \bbeta_{10}) \stackrel{{\cal L}} \longrightarrow {\cal N}({\bf 0}, {\cal I}_{10}^{-1}[{\cal D}_{10}^T{\cal W}_{10} \bSigma_{0}^{(1)}{\cal W}_{10}{\cal D}_{10}]{\cal I}_{10}^{-1}),$$ under $P_{\bbeta_0 \pi}$, showing that the penalized pseudo-QIF estimator is as efficient as the oracle estimator that assumes the true model, with $\bbeta_{20}={\bf 0}$ is known.
}
\end{remark}

By Theorem \ref{normality}, an estimator of the asymptotic variance of $\tilde{\bbeta}_{1n}$ is given by 
$$n^{-1} \left[\frac{\partial^2 n^{-1}Q_n}{\partial \bbeta_1 \partial \bbeta_1^T} \left\{\begin{pmatrix}
\tilde{{\bbeta}}_{1n} \\
{\bf 0}
\end{pmatrix} \right\}  +{\bf B}_n\right]^{-1} \hat{{\cal V}}_0^{(1)} \left[\frac{\partial^2 n^{-1}Q_n}{\partial \bbeta_1 \partial \bbeta_1^T} \left\{\begin{pmatrix}
\tilde{{\bbeta}}_{1n} \\
{\bf 0}
\end{pmatrix} \right\} +{\bf B}_n \right]^{-1},$$
where $\hat{{\cal V}}_0^{(1)} =  4{\cal D}_{1n}(\tilde{\bbeta}_{1n})^T\left[{\bf C}_{1n}(\tilde{\bbeta}_{1n})\right]^{-1} \hat{\bSigma}_{0}^{(1)}\left[{\bf C}_{1n}(\tilde{\bbeta}_{1n})\right]^{-1}{\cal D}_{1n}(\tilde{\bbeta}_{1n}),$ with $\hat{\bSigma}_{0}^{(1)} = \hat{\bSigma}_d +f \hat{\bSigma}_{(\tilde{\bbeta}_{1n}^T, {\bf 0})^T}.$ In addition, ${\cal D}_{1n}(\tilde{\bbeta}_{1n})$ denotes the $Ld_1 \times d_1$ matrix obtained from ${\cal D}_n\left\{\begin{pmatrix}
\tilde{\bbeta}_{1n} \\
{\bf 0}
\end{pmatrix}\right\}$ by selecting the first $d_1$ columns and rows numbered by $\eta d+1, \ldots, \eta d+d_1$ and ${\bf C}_{1n}(\tilde{\bbeta}_{1n})$ is the $Ld_1 \times Ld_1$ matrix obtained from ${\bf C}_n\left\{\begin{pmatrix}
\tilde{\bbeta}_{1n} \\
{\bf 0}
\end{pmatrix}\right\}$ by selecting the rows and columns which are numbered as $\eta d+1, \ldots, \eta d+d_1$, where $\eta=0, \ldots, (L-1)$. 

In complex sample surveys the asymptotic variance of an estimator may have a complicated form and resampling methods have to be used. \cite{wang-wang-wang14} used the estimating function bootstrap method for variance estimation under the penalized pseudo-GEE method for variable selection. A bootstrap method can be employed to obtain a variance estimator of $\tilde{\bbeta}_{1n}$ by generating bootstrap weights using the Rao-Wu rescaling method (see \cite{rao-wu88}, \cite{rao-wu-yue92}) and taking a one-step bootstrap. Then, as in \eqref{NRupdate}, for each set of bootstrap weights, the value of the estimator of the non-zero components can be updated as follows
\begin{equation}
\tilde{\bbeta}_{1n}^{(b)} = \tilde{\bbeta}_{1n} - \left[ \frac{\partial^2 Q_n^{(b)}(\tilde{\bbeta}_n)}{\partial \bbeta_1 \partial \bbeta_1^T} + n \boldsymbol{\Gamma}(\tilde{\bbeta}_n)\right]^{-1} \left[\frac{\partial Q_n^{(b)}(\tilde{\bbeta}_n)}{\partial \bbeta_1}+ n\boldsymbol{\Gamma}(\tilde{\bbeta}_{n})(\tilde{\bbeta}_{1n})\right], \label{boot-penal-estim}
\end{equation}
where $Q_n^{(b)}(\tilde{\bbeta}_{1n})$ is the bootstrap weighted QIF, calculated using the $b$-th set of bootstrap weights. Consequently, a bootstrap variance estimator is given by 
\begin{equation}
\hat{{\bf V}}^B_{\tilde{\bbeta}_{1n}} :=\frac{1}{B}\sum_{b=1}^B(\tilde{\bbeta}_{1n}^{(b)}-\tilde{\bbeta}_{1n})(\tilde{\bbeta}_{1n}^{(b)}-\tilde{\bbeta}_{1n})^T.\label{boot-var}
\end{equation}

\section{Numerical results}
\label{simulations}

We generate a finite population with correlated binary responses from the following marginal logistic model
$$ {\rm logit} \mu_{ij}(\bbeta_0) = {\bf x}_{ij}^T\bbeta_0, \ {x}_{ij}^{k} \stackrel{indep}\sim {\cal U}(0, 0.8), \ i=1, \ldots N, \ j=1, \ldots, m, \ k=1, \ldots d,$$
 $N=30000,$ $m=5,$ $d=10$, $\bbeta_0=(0.8, -0.7, -0.6, 0 , 0 , 0, 0, 0, 0, 0)^T$, choosing the exchangeable correlation matrix with parameter $\alpha=0.4$ as the true correlation. In this case, the basis matrices are ${M}_1={\bf I}_5$ and ${\bf M}_2=\begin{bmatrix}
0 & 1 & 1 & 1 & 1 \\
1 & 0 & 1 & 1 & 1 \\
1 & 1 & 0 & 1 & 1 \\
1 & 1 & 1 & 0 & 1 \\
1 & 1 & 1 & 1 & 0 \\
\end{bmatrix}$ and $d_1=3$. 

From each generated finite population we obtain a sample of $n$ clusters, using informative sampling: clusters are selected with probability proportional to the size measures $z_i=\sum_{j=1}^5 y_{ij}+1$, with replacement. We consider sample sizes of $n=300$ and $n=500$ and repeat the procedure of generating a finite population and then selecting the sample $H=500$ times. 

For each sample, we apply each the following methods: unweighted QIF with SCAD penalty (UNWGT), weighted QIF with SCAD penalty (PQIF), weighted GEE with SCAD penalty (PGEE) and ORACLE, which is the weighted QIF under the true model (with three nonzero coefficients and seven zero coefficients). 

To evaluate the performance of the proposed method, two working correlations: exchangeable (EX) and fist-order autoregressive (AR1) are considered for each procedure and assessed in Table \ref{sim-table}, as follows. The columns labeled $(C)$, $(O)$ and $(U)$ give, respectively, the percentage of times the true model (only first three components of the estimator are non-zero) is selected, the percentage of times the variables are over-selected (more than the first three components of the estimator are non-zero ) and the percentage of times the variables are under-selected (at least one of the first three components of the estimator is zero). The results in Table \ref{sim-table} show that the unweighted PQIF is not capable of variable selection and only less than 10\% of the time the true model selected. Furthermore, PQIF and PGEE both perform well in terms of selecting the true model (C), with PQIF leading to slightly larger values.

In addition, we compute the Monte Carlo average MSE of the estimators, reported in Table \ref{sim-table} as MSE
$$ MSE=\frac{1}{H} \sum_{h=1}^{H} (\tilde{\bbeta}_n^{(h)} - \bbeta_0)^T(\tilde{\bbeta}_n^{(h)} - \bbeta_0),$$
where $\tilde{\bbeta}_n^{(h)}$ is the survey-weighted penalized estimator, calculated at iteration $h,$ with $h=1, \ldots, H$ and $H=500$. Table \ref{sim-table} shows that UNWGT leads to larger MSE, while PQIF exhibits smaller MSE than PGEE. As expected, ORACLE performs the best in terms of correctly selecting the true model and MSE. 

\begin{table}[htbp] 
\caption{Correlated binary responses in surveys: comparisons of the unweighted QIF, penalized pseudo-QIF, penalized pseudo-GEE and oracle QIF under exchangeable and AR(1) working correlation matrices. The columns labeled $(C)$ indicate the percentage of times the true model is selected (i.e. only first three components of the estimator are non-zero). The columns labeled $(O)$ indicate the percentage of times the variables are over-selected (i.e. more than the first three components of the estimator are non-zero). The columns labeled $(U)$ give the percentage of times the variables are under-selected (i.e. at least one of the first three components of the estimator is zero).}
\small{\begin{center}
	{\begin{tabular}[c]{l l r r r r r r r r r r r r }
	\toprule
	Sample      & Method &   \multicolumn{4}{c}{Exchangeable} & \multicolumn{4}{c}{AR(1)}   \\
	                      \cmidrule(lr){3-6}         \cmidrule(lr){7-10}
		size			&        &    (C) &   (O)  &  (U)   & MSE   &   (C) &   (O)  &  (U)   & MSE  \\		
	\midrule
	$n=300$     &  UNWGT & 5.6 & 34.2 & 60.2 & 0.538 & 2.4 & 16.0 & 71.6 & 0.610 \\	
              &  PQIF  & 80.8& 11.2 & 8.0  & 0.129 & 75.0& 11.4 & 13.6 & 0.163 \\
							&  PGEE  & 73.0& 10.6 & 16.4 & 0.172 & 72.2& 10.4 & 17.4 & 0.181 \\
              & ORACLE &100.0& 0.0  & 0.0  & 0.072 &100.0& 0.0  & 0.0  & 0.080 \\
							\\
	$n=500$     & UNWGT  & 7.0 & 66.6 & 26.4 & 0.384 & 5.4 & 59.6 & 35.0 & 0.435 \\
              & PQIF   & 91.2& 8.6  &  0.2 & 0.049 &87.8 & 11.0 & 1.2  & 0.065 \\
							&  PGEE  & 88.2& 10.0 & 1.8  & 0.054 & 87.0& 10.2 & 1.8  & 0.062 \\
              & ORACLE &100.0& 0.0  &  0.0 & 0.039 &100.0& 0.0  & 0.0  & 0.046 \\		
\midrule
				\end{tabular}
				}
\end{center}}
\label{sim-table}
\end{table}

In Table \ref{bias-table} we report the absolute relative bias of the non-zero coefficients, calculated as  
$ \displaystyle{ARB(\tilde{\beta}_{nk}) = \frac{\left|H^{-1} \sum_{h=1}^H \tilde{\beta}_{nk}^{(h)} - \beta_{0k} \right|}{|\beta_{0k}|} \times 100, \ k=1,2,3}$, where $\tilde{\beta}_{nk}^{(h)}$ is the $k$-th component of $\tilde{\bbeta}_{n}^{(h)}$ and $\beta_{0k}$ is the $k$-th entry of $\bbeta_0$. The unweighted PQIF yields biased results, whereas the relative biases of PQIF estimates are less that those of PGEE. In addition, in contrast to the unweighted methods, the values the relative bias decrease in case of the weighted methods.

\begin{table}[htbp] 
\caption{Correlated binary responses in surveys: comparisons of percent absolute relative bias (ARB) of regression coefficients, obtained from the unweighted QIF, penalized pseudo-QIF, penalized pseudo-GEE and oracle QIF, under exchangeable and AR(1) working correlation matrices.}
\small{\begin{center}
	{\begin{tabular}[c]{l l r r r r r r r r r r }
	\toprule
	Sample      & Method &   \multicolumn{3}{c}{Exchangeable} & \multicolumn{3}{c}{AR(1)}   \\
	                      \cmidrule(lr){3-5}         \cmidrule(lr){6-8}
		size			&        & $\tilde{\beta}_{n1}$ & $\tilde{\beta}_{n2}$ & $\tilde{\beta}_{n3}$ &  $\tilde{\beta}_{n1}$ & $\tilde{\beta}_{n2}$ & $\tilde{\beta}_{n3}$   \\		
	\midrule
	$n=300$     &  UNWGT &28.8  &30   & 55   &  28.8 & 37.1 & 61.7  \\	
              &  PQIF  & 1.3  & 1.4 &  6.7 &   2.5 &  0   &  8.3  \\
							&  PGEE  & 7.5  & 1.4 & 13.3 &   6.3 &  1.4 & 13.3  \\
              & ORACLE & 1.3  & 1.4 &  3.3 &   0   &  0   &  3.3  \\
							\\
	$n=500$     & UNWGT  & 28.4 & 22.9& 41.4 &  28.6 &26.6  & 46.5  \\
              & PQIF   & 1.7  & 2.3 & 4.0  &   1.4 & 2.0  & 3.9  \\
							&  PGEE  & 3.9  & 1.2 & 5.8  &   3.5 & 0.9  & 4.5   \\
              & ORACLE & 1.8  & 2.4 & 3.9  &   1.5 & 2.0  & 3.2   \\	
\midrule								
				\end{tabular}
				}
\end{center}}
\label{bias-table}
\end{table}

Next, we evaluate the performance of the bootstrap variance estimator given by \eqref{boot-var}. To obtain the bootstrap standard error, we draw 500 bootstrap samples of size $n-1$, each selected with replacement and equal probabilities from the $n$ sampled units.
For each $b$-th bootstrap sample, the bootstrap weights are calculated using the rescaling formula
$$w_i^{(b)} = w_i\frac{n}{n-1} t_i^{(b)},$$
where $t_i^{(b)}$ is the number of repetitions of unit $i$ in the $b$-th bootstrap sample. 

As in \cite{fan-li01}, to illustrate the performance of the proposed variance estimator, we first compute $SD$ (for each $k=1,2,3$), defined as the ratio between the median absolute deviation of the estimator $\tilde{\beta}_{nk}$
$$\mbox{median} \left\{\left|\tilde{\beta}_{nk}^{(h)} - \mbox{median}\{\tilde{\beta}_{nk}^{(h)},  \ h=1, \ldots, H \} \right|, \ h=1, \ldots, H \right\},$$
and 0.6745. Due to the normality of the limiting distribution of $\tilde{\beta}_{nk}$, this value is an estimate of its true standard error and in Table \ref{var-estim}, we compare it with the one obtained from the bootstrap method, denoted as $SD_m$. The latter is obtained as the median of the $H$ bootstrap estimated standard deviations
$$\mbox{median} \left\{\sqrt{\hat{{\bf v}}^{B,h}_{\tilde{\beta}_{nk}}}, \ h=1, \ldots, H \right\}, $$
where, at each iteration $h$, $\displaystyle{\hat{{\bf v}}^{B,h}_{\tilde{\beta}_{nk}} =\frac{1}{B}\sum_{b=1}^B(\tilde{\beta}_{nk}^{(b, h)}-\tilde{\beta}_{nk}^{(h)})^2}$ is the bootstrap variance estimate of $\tilde{\beta}_{nk}^{(h)}$. Here, $\tilde{\beta}_{nk}^{(b, h)}$ denotes the $k$-th component of $\tilde{\bbeta}_{1n}^{(b,h)}$ whose form is given in \eqref{boot-penal-estim}, $h=1, \ldots, H$. Furthermore, for each $k=1,2,3$, the median absolute deviation error of the 500 bootstrap estimated standard errors, denoted as $SD_{mad}$ is evaluated as the ratio between 
$$\mbox{median} \left\{\left|\sqrt{\hat{{\bf v}}^{B,h}_{\tilde{\beta}_{nk}}} - SD_k \right|, \ h=1, \ldots, H \right\},$$
and 0.6745, where $SD_k$ is the estimate of the true standard error of $\tilde{\beta}_{nk}$. Results in Table \ref{var-estim} compare estimates obtained under the two working correlations and different sample sizes. They show that in each case, the one-step bootstrap method for variance estimation performs well in tracking the value of the true standard error of $\tilde{\beta}_{nk}$, $k=1,2,3$. 

\begin{table}[htbp] 
\caption{Correlated binary responses in surveys: accuracy of the bootstrap approximation of the variance of the penalized pseudo-QIF estimator, under  exchangeable and AR(1) working correlation matrices.}
\small{\begin{center}
	{\begin{tabular}[c]{l  r r  r r r r  }
	\toprule
	              &  \multicolumn{2}{c}{$\tilde{\beta}_{n1}$}& \multicolumn{2}{c}{$\tilde{\beta}_{n2}$} & \multicolumn{2}{c}{$\tilde{\beta}_{n3}$} \\
	Sample				&  SD & $SD^B_m$ $(SD^B_{mad})$ & SD & $SD^B_m$ $(SD^B_{mad})$ & SD & $SD^B_m$ $(SD^B_{mad})$ \\ 
  						 \cmidrule(lr){2-7} 
							\\
	size 						& \multicolumn{6}{c}{Exchangeable}  \\
	\midrule
	$n=300$     &  0.172 & 0.163 (0.016) & 0.152 & 0.161 (0.015) & 0.142 & 0.156 (0.015) \\	
 							\\
	$n=500$     &  0.112 & 0.127 (0.011) & 0.125 & 0.125 (0.010) & 0.099 & 0.123 (0.010)  \\
 \midrule			\\				
              &   \multicolumn{6}{c}{AR(1)}  \\
	\midrule
	$n=300$     &  0.181 & 0.172 (0.018) & 0.170 & 0.170 (0.016) & 0.162 & 0.165 (0.016) \\	
 							\\
	$n=500$     &  0.124 & 0.134 (0.011) & 0.137 & 0.133 (0.011) & 0.125 & 0.131 (0.011)  \\
	\midrule
				\end{tabular}
				}
\end{center}}
\label{var-estim}
\end{table}

\vspace{3mm}

\bibliographystyle{plainnat}
	
	\end{document}